\documentclass[12pt]{amsart}
\usepackage{graphicx}
\usepackage{hyperref}
\usepackage [latin1]{inputenc}

\newtheorem{theorem}{Theorem}
 
\theoremstyle{definition}

\newtheorem{example}{Example}
\newtheorem{remark}{Remark}
\def\r{\mathbb R}

\title{What is the Shape of a Cupola?}
\author{Rafael L\'opez}
\address{Departamento de Geometr\'{\i}a y Topolog\'{\i}a\\ Universidad de Granada. Granada, Spain}
\email{rcamino@ugr.es}

\begin{document}

\begin{abstract}  This article examines the shape of a surface obtained by a hanging flexible, inelastic material with prescribed area and boundary curve. The shape of this surface, after being turned upside down, is a model for cupolas (or domes) under the simple hypothesis of compression.   Investigating the rotational examples, we provide and illustrate a novel design for a roof   which has the extraordinary property that its shape, although natural, is modeled by a surface of revolution whose axis of rotation is horizontal.

\end{abstract}
 
\date{}
\maketitle

\section{Introduction.}

Historically, the shape of a cupola (or dome) has been of enduring interest. The Greek's use of columns and the Roman's use of arches as a basic element in construction enabled architects  to build ever larger walls and pillars, increasing the relevance of the cupola as the crowning element of the entire edifice. The use of flying buttresses to distribute loads and tensions in walls over a large area transformed the low windowless Romanesque churches into the tall, slender Gothic cathedrals that embellish the cities of Europe.

\begin{figure}[hbtp]
\begin{center}
\includegraphics[width=.3\textwidth]{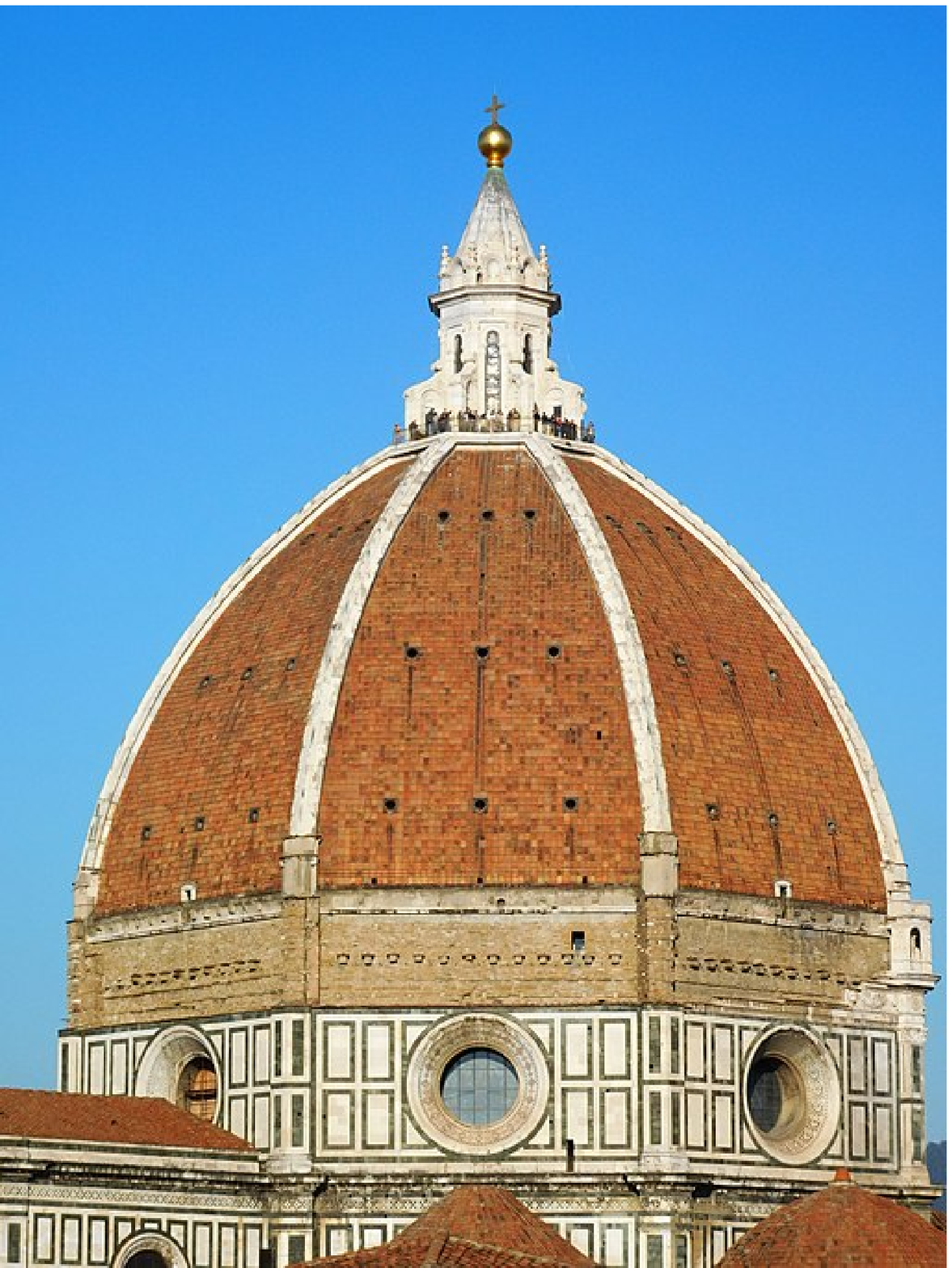}\qquad \includegraphics[width=.5\textwidth]{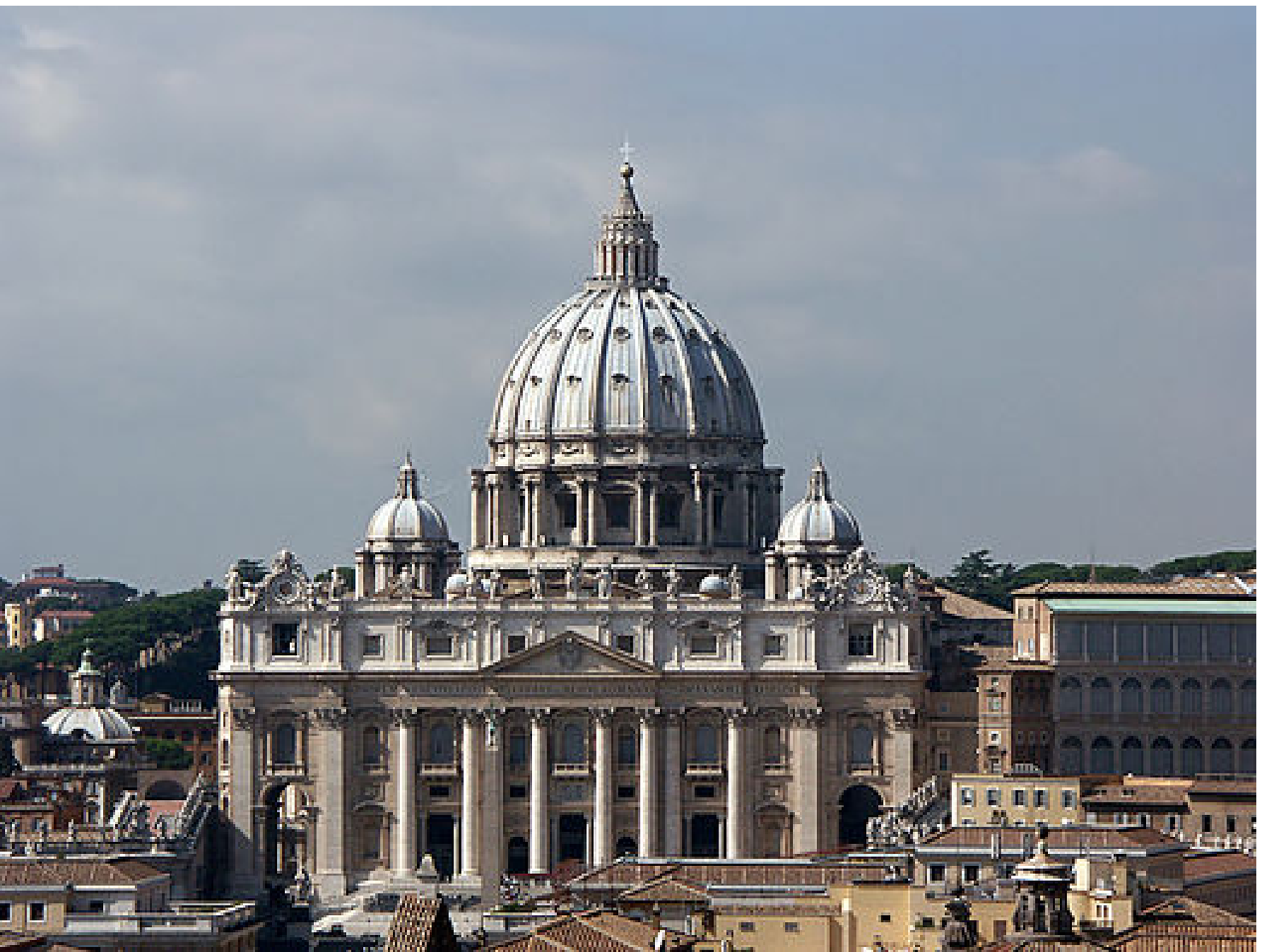}
\end{center}
\caption{Domes of Santa Maria del Fiore in Florence (left) and  the Basilica Papale di San Pietro in Vaticano (right). The first image is licensed under the Creative Commons Attribution-Share Alike 4.0 International license at commons.wikimedia.org/wiki/File:Florence duomo fc10.jpg. The second image is licensed under the public domain at  commons.wikimedia.org/wiki/File:Petersdom von Engelsburg gesehen.jpg.}\label{fig1}
\end{figure}

The construction of cupolas involves an intricate interplay of artistic and structural issues requiring  the architect to specify a variety of variables such as the choice of materials and the desired stylistic effect. The essential engineering problem to be solved is to build a large structurally stable, aesthetically appealing roof that rises over a large, empty space. In order to achieve this, architects certainly required the support of the sciences.  In the 15th and 16th centuries, the Renaissance was a period of scientific and artistic development propitious for the building of domes. The examples of the Florentine Cathedral of Santa Maria del Fiore by Filippo Brunelleschi and the Vatican's Basilica Papale di San Pietro of Michelangelo (Figure \ref{fig1}) demonstrate the resulting  triumphant achievements that fascinate us up to the present day.

Owing to the issues outlined above, it is not clear how one should go about  formulating the problem of finding the optimal shape for a cupola. Here we take a mathematician's perspective.  A first thought that comes to mind is that the cupola is sustained along its boundary by its own weight.
As a first approximation, we imagine a bounded, massive, homogenous  piece of cloth whose boundary is represented by a fixed prescribed curve. Supported by this curve, the cloth evolves under the force of gravity to a static equilibrium.  The reason the surface of the cloth is closely related to that of a  cupola is that a surface suspended solely by its own weight experiences only tensional forces tangent to its interior. When this surface is inverted, it produces the optimal shape of a cupola. The inversion transforms the tensional force into a force of compression. In our context of cupolas, we can rephrase the words  of Robert Hooke about the shape of an arch by saying that as hangs the flexible surface so inverted will stand the rigid cupola.\footnote{Actually Hooke considered the problem of the hanging cable writing an anagram in Latin that deciphers to ``{\it ut pendet continuum flexile, sic stabit contiguum rigidum inversum}'' which translates as  ``as hangs the flexible line, so but inverted will stand the rigid arch'' (\cite[p. 31]{ho}).} The inverted surface satisfies the same equation of equilibrium as the original surface so the only question to be answered is:

\medskip
\noindent {\it What is the shape of a flexible hanging surface of uniform mass acted upon solely by gravity?}
\medskip
  
\begin{figure}[hbtp]
\begin{center}
\includegraphics[width=.34\textwidth]{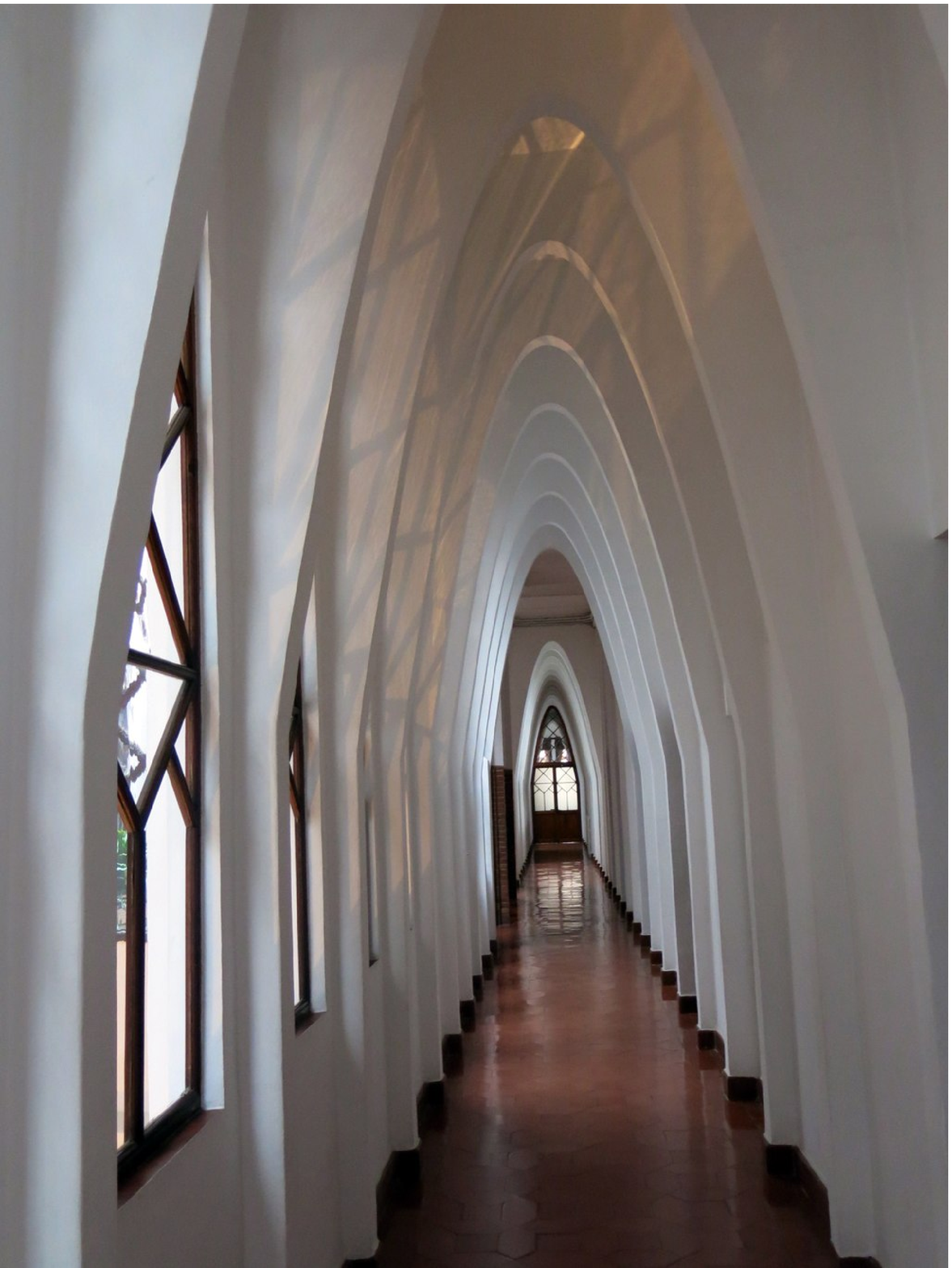}\qquad \includegraphics[width=.3\textwidth]{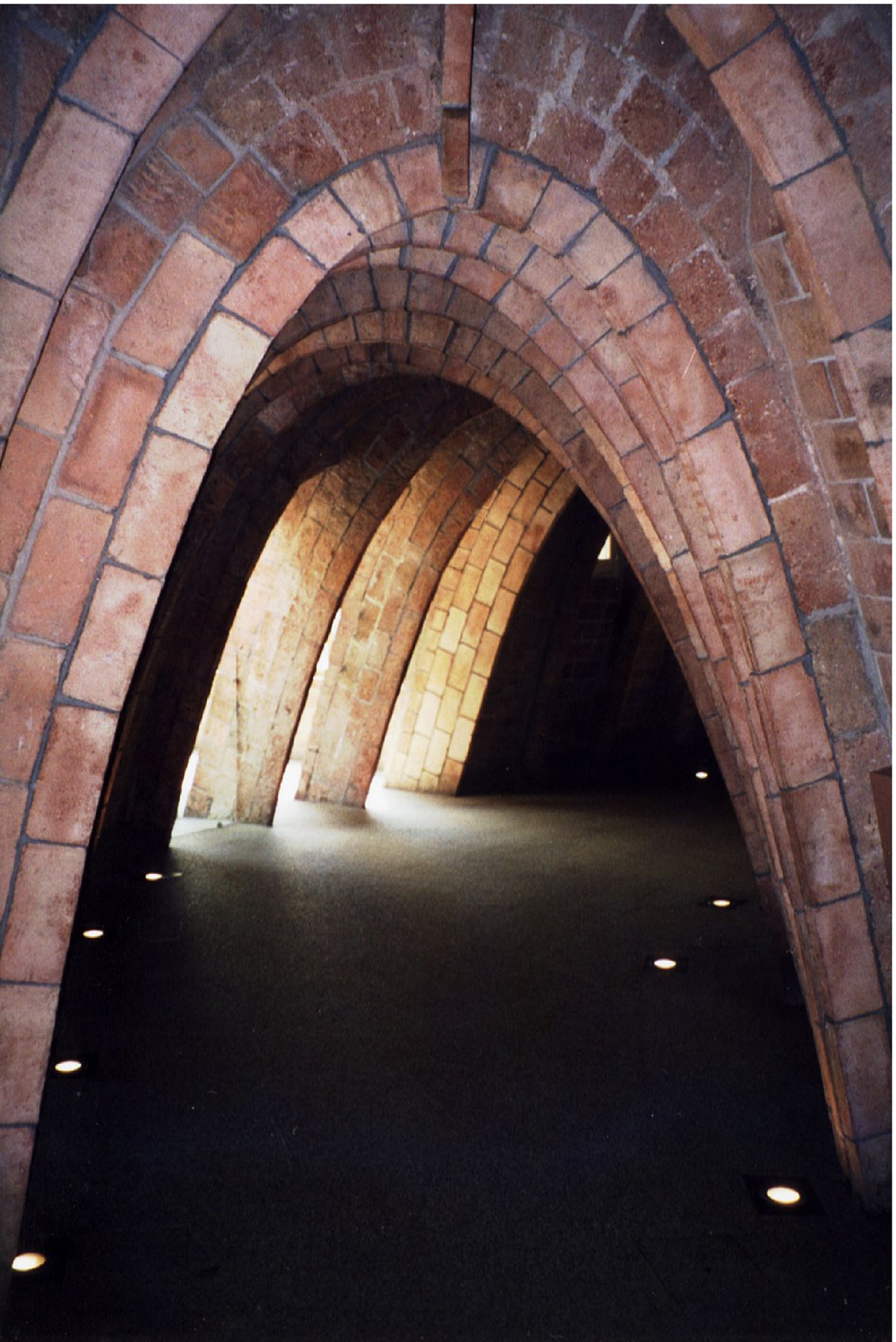}
\end{center}
\caption{Left: Corridor in the Colegio Teresiano, Barcelona. Right: loft in La Pedrera, Barcelona. The first image   is licensed under the Creative Commons Attribution-Share Alike 3.0 Unported license at commons.wikimedia.org/wiki/File:049 Col·legi de les Teresianes, arcs parabòlics.JPG. The second image is licensed under the GNU Free Documentation License at commons.wikimedia.org/wiki/File:LaPedreraParabola.jpg.}\label{fig2}
\end{figure}

  As is often the case, some insight can be achieved by considering the one dimensional analog of the problem stated above, which is to determine  the optimal shape of a hanging cable.  The answer, as is well known, is a catenary curve given by the simple expression $y(x)=a^{-1}\cosh(ax)$, for a positive constant $a$. The optimal shape of arches has also attracted the interest of mathematicians, where catenaries and parabolas have competed for this role; see, for example the beautiful discussion of R. Osserman on the shape of the Gateway Arch in Saint Louis, Missouri (\cite{os}).  The  renowned Spanish architect Antonio Gaud\'{\i} (1852-1926), who included many beautiful catenary shaped corridors, was an avid enthusiast of this shape (see Figure \ref{fig2}).

The list of mathematicians who have investigated the shape of surfaces hanging under their own weight includes the names of Beltrami, Germain, Jellet, Lagrange and Poisson (\cite{be,ge,je,la,po,vo}). Surely, it was Cisa de Gresy who stated  most lucidly \cite[p. 260]{ci}:

\begin{quote}
  ``Si on suppose, par example, une surface en \'equilibre, sollicit\'ee uniquement par la gravit\'e, et suspendue \`a la circonf\'erence d'un cercle fix\'e horizontalement, it est clair que les \'el\'ements de cette surface n'\'eprouveront qu'une simple tension dans le sens des m\'eridiens ou de la courbe g\'en\'eratrice.'' [If we suppose, for example, that a surface is subjected  only to the force of gravity and it is suspended from a circular perimeter, it is clear that the elements of this surface will only exert a simple tension in all directions of the meridians or the generating curve.]
\end{quote}

However, it is possible that there is not a minimum of the height of  the center of gravity for surfaces with prescribed area and  boundary curve. We present an example which is a   slightly simplified 
version of the one given by Nitsche in \cite{ni}. 

\begin{example}\label{ex1} Let $\Gamma$ be the circle in the plane $z=0$ of radius $1$ and centered at the origin. For $0<R<1$, let   $\Omega_R$ be the annulus $\{(x,y,0)\in \r^3: R^2\leq x^2+y^2\leq 1\}$. Consider the surface $S_R$ formed by $\Omega_R$ together with the  cone $C_R$ underneath $\Omega$ with boundary $C_R=\{(x,y,0)\in\r^3: x^2+y^2=R^2\}$ and height $h=\sqrt{2+1/R^2}$. We can parametrize $S_R$ in polar coordinates $(r,\theta)$ by 
$$u(r)=\left\{\begin{array}{lll} -h\frac{R-r}{R},& & 0\leq r<R\\
0, & & R\leq r\leq 1.\end{array}\right.$$
See Figure \ref{figcone}. The boundary of $S_R$ is the circle $\Gamma$ and with these choices of $R$ and $h$, the area of $S_R$ is constantly $2\pi$ independently of $R$ (the value $2\pi$ is only for convenience; any area greater than $\pi$ can be taken). It is well known that the center of gravity of the (hollow) cone of height $h$ is $h/3$ from the base.   So the center
of gravity of $S_R$ is at height 
\[
	\frac{-\frac{h}{3}\cdot\mbox{area}(C_R)}{\mbox{area}(S_R)}=\frac{-\frac{h}{3}\pi R\sqrt{h^2+R^2}}{2\pi}=-\frac{1+R^2}{6}\sqrt{2+\frac{1}{R^2}}.
\]
\begin{figure}[h]
\begin{center} 
\includegraphics[width=.7\textwidth]{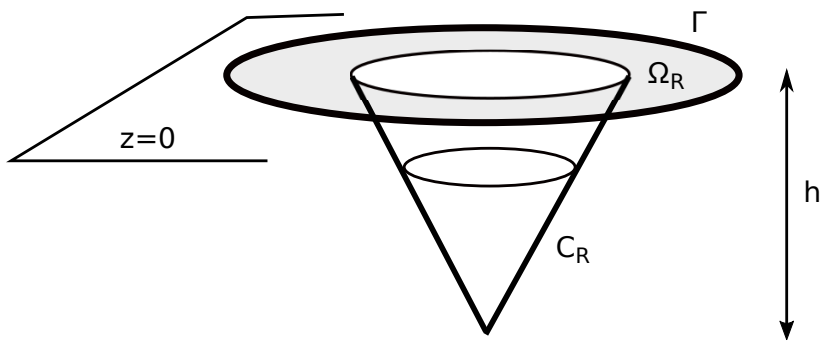}%
\end{center}
\caption{The problem of minimizing the height of the center of gravity has no solution when the prescribed area is $2\pi$ and the boundary curve is the circle $\Gamma$. See Example \ref{ex1}.} \label{figcone}
\end{figure}
In   particular, the center of gravity can be made as low as one likes, by taking $R$ sufficiently small.  This example can
obviously also be made smooth by small modifications. 
\end{example}

This example contrasts with the one dimensional version of the problem because, as was proved by Jacob Bernouilli, the catenary has the property that its center of gravity is lower than   that of any curve of equal length, and with the same fixed endpoints.

\begin{figure}[h]
\begin{center}
\includegraphics[width=.95\textwidth]{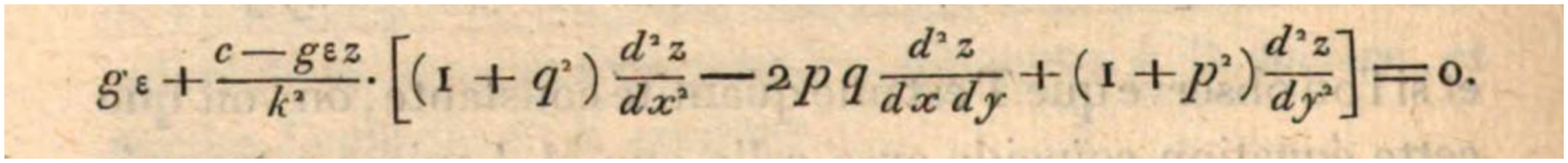} 
\end{center}
\caption{The equation deduced by Poisson that satisfies a surface $z=(x,y)$ acted upon solely by gravity. } \label{figpo}
\end{figure}

As is usual in optimization problems, and in light of the above example, we approach the
problem by requiring something less than an absolute minimum for the height of
the center of gravity. Indeed, when a (flexible, inelastic) material hangs under its
weight, the surface that is formed is a local extremum for the height of the center
of gravity, in the space of smooth surfaces with given area and boundary.  Therefore techniques from the Calculus of Variations are key for deriving the differential equation of the surface. In relation to this,  Joseph-Louis Lagrange states in his \textit{M\'ecanique Analytique}:
  \begin{quote}
``[...] on verra par l'uniformit\'e et la rapidit\'e des solutions combien ces m\'ethodes sont sup\'erieures \`a celles que l'on avait employ\'ees jusqu'ici dans la Statique''. [[...] one will see by the consistency and speed of solvability, how  these methods are greater than to those that have been employed until now in Statics.]  See \cite[p. 113]{la}.
\end{quote}
 But it was Sim\'eon Denis Poisson who definitively found the equilibrium equation for the surface, improving the assumptions and calculations of Lagrange. What's more, Joseph Bertrand, who edited the collected works of Lagrange, added a footnote:
\begin{quote}    
`Cette mani\`ere d'\'evaluer l'ensemble des forces que d\'eveloppe l'\'elasticit\'e sur un point n'est pas suffisamment justifi\'ee [...] Nous pouvons m\^{e}me ajouter que cela n'est pas exact. Poisson en a fait la remarque dans le {\it M\'emoires de l'Institut} pour l'ann\'ee 1812''. [This way of evaluating the collection of forces, which develop the elasticity at a point, is not sufficiently justified [...]. We may even add that it is not exact. Poisson made this observation in {\it M\'emoires de l'Institut} in 1812]. See   \cite[p. 158]{la}.
\end{quote}
  Indeed, Poisson   considered a much more general problem of a surface under different forces and tensions. As a particular case, he derived the correct equation of the surface stretched by its weight, which we will see in the next section. So, assuming only the effect of the weight, he asserts:

 \begin{quote} ``Consid\'erons enfin la surface pesante, et prenons l'axe des $z$ vertical et dirig\'e dans le sens de la pesanteur''. [Let us finally consider the heavy surface. We take the vertical axis pointing along the direction of the gravitational field.]  See \cite[p. 185]{po}.
 \end{quote}
Then he successfully derived the equation for a nonparametric surface $z=z(x,y)$ (see Figure \ref{figpo}), where $k^2=1+p^2+q^2$, $p=z_x$, $q=z_y$, $g$ is the gravitational acceleration and $\epsilon$ is the density of the surface.
  Finally, he writes: 
 \begin{quote} ``Cette \'equation d'\'equilibre de la surface pesante et \'egalement \'epaisse, doit comprendre l'\'equation ordinaire de la {\it cha\^{i}nette}, qui s'en d\'eduit, en effect, en y supposant $z$ ind\'ependante de l'une des deux variables $x$ ou $y$, de $y$, par exemple''. [This equilibrium equation of the heavy surface with uniform thickness must include the known equation of the {\it suspended chain}, which is deduced from it by assuming that $z$ is independent of one of the two variables $x$ or $y$, say of $y$]. See  \cite[p.186]{po}. 
 \end{quote} 
All the aforementioned works were apparently nearly forgotten until the 1980's, when there was an explosion of interest in the evolution of surfaces by functions of their mean curvature.  There is also the issue of the elasticity of the materials used in the construction. As the reader can well imagine, a dome's actual material is not nearly so flexible as the cloth example discussed above.  See a historical approach   in  \cite{to}.  Here, we would like to take note of the paper \cite{di2} by Ulrich Dierkes,   which was surely motivated by the work of the German architect Frei Otto (\cite{ot}). Later, the problem was revisited by Bemelmans, B\"{o}hme, Dierkes, Hildebrandt, and Huisken    in  their works (\cite{bu,bht,di1,di2,dh}).

The literature in architecture on the shape of cupolas is extensive and cannot be catalogued  here. We refer only to \cite{du,he,op,pa,pot}. Since we lack expertise in the fields of architecture and engineering, we have approached the problem from the perspective of differential geometry, although we have avoided its technical concepts, such as shape operator, principal curvatures, and second fundamental form, to maintain accessibility for a larger readership.

The surfaces we discuss below will be graphs, surfaces of revolution, or cylindrical surfaces whose parameterizations are simple. In Section \ref{sec2} we will employ the calculus of variations to derive the equation that a function  $z=u(x,y)$  must satisfy for its graph to define a surface whose shape is determined only by its own weight. These surfaces are called {\it singular minimal surfaces}.   We will see that the boundary of the surface imposes geometric restrictions to the shape of the entire surface and this question will be briefly discussed. In Section \ref{sec3} we focus on singular minimal surfaces that are surfaces of revolution, thinking of the shape of cupolas. Finally, in   Section \ref{sec4} we will present a new roof design modeled by a singular minimal surface. The novelty is that the roof is a surface of revolution but its rotation axis is horizontal,   which is contrary to our common sense.

\section{Singular minimal surfaces.}\label{sec2}

Consider $(x,y,z)$ the canonical coordinates of the three-dimensional Euclidean space $\r^3$  where $z$ indicates the vertical direction. Let $\Gamma$ be a closed curve and $A>0$ a fixed positive number. We wish to determine the differential equation that governs a surface $S$ spanning $\Gamma$ with area $A$ which is suspended from $\Gamma$ by its weight. Suppose that  $S$  is made of a  flexible, incompressible material of uniform density $\sigma$ per unit area.  In order to simplify the arguments, we restrict our attention to   surfaces  given by the graph of a smooth function $z=u(x,y)$ defined on $\Omega$, a bounded planar domain with smooth boundary $\partial\Omega$.  The weight per unit area of $S$  is $\sigma\sqrt{1+u_x^2+u_y^2}$, where the subscripts indicate the derivatives with respect to the corresponding variables.  Under the effect of the weight, the surface $S$ attains a  point of equilibrium when the height of its center of gravity is a local extremum. Assume that  the gravitational potential at one point $(x,y,z)$ is simply the distance $z$ to the $xy$-plane. In particular, all our geometric objects (curves and surfaces) lie over the plane of equation $z=0$.   Let us also observe that the problem is invariant under translations in any horizontal direction. The height of the center of gravity  is 
$$\frac{1}{A}\int_{\Omega}\sigma\ u\sqrt{1+u_x^2+u_y^2}\, dxdy.$$
 The minimization is understood to be in the class of smooth functions $u$ with prescribed boundary $u=\varphi>0$, where the graph of    $\varphi\colon \partial\Omega\to\r$ is just the boundary curve $\Gamma$.  We can assume that $A$ and $\sigma$ take the value $1$. 
 
We now consider simple arguments of calculus of variations and make an infinitesimal change in the surface $z=u(x,y)$ given by 
  $u(x,y)+t h(x,y)$, $t\in\r$ and $h\colon\Omega\to\r$ a smooth function vanishing on $\partial\Omega$.  Adopting a Lagrange multiplier for the constraint on the area of the surface, define the functional 
\begin{equation}\label{jj}
J(u)=\int_\Omega u\sqrt{1+|\nabla u|^2}\, dxdy+\lambda \int_\Omega \sqrt{1+|\nabla u|^2}\, dxdy,
\end{equation}
where $\nabla u=(u_x,u_y)$ stands for the gradient of $u$ and $\lambda\in\r$.  The domain of $J$ is the set of all smooth functions $u$ defined on $\Omega$ with boundary condition $u=\varphi$ along $\partial\Omega$ and fixed surface area equal to $1$. The class $\mathcal{X}$ of admissible variations   is formed by the smooth functions $h\colon\Omega\to\r$   which vanish on the boundary of $\Omega$, $h=0$ along $\partial\Omega$. Thus an   extremal $u$ of $J$ implies    
$$\frac{d}{dt}{\Big|}_{t=0} J(u+th)=0$$
for all $h\in \mathcal{X}$.  Set the Lagrangian $L(x,y,u,p,q)=(u+\lambda)\sqrt{1+p^2+q^2}$, the integrand in  \eqref{jj}, with $p=u_x$ and $q=u_y$,  and let 
$$F=\left(\frac{\partial L}{\partial p},\frac{\partial L}{\partial q}\right)=(u+\lambda)\left(\frac{u_x}{\sqrt{1+|\nabla u|^2}},\frac{u_y}{\sqrt{1+|\nabla u|^2}}\right).$$
Using  $\mbox{div} (h\cdot F)=\langle\nabla h,F\rangle+h\cdot \mbox{div}F$, where $\langle\cdot,\cdot\rangle$ denotes the usual scalar product, we have
\begin{equation*}
\begin{split}
\frac{d}{dt}{\Big|}_{t=0} J(u+th)&=\int_\Omega\left(\frac{\partial L}{\partial u}h+h_x\frac{\partial L}{\partial p}+h_y\frac{\partial L}{\partial q}\right)\, dxdy\\
&=\int_\Omega\left(\frac{\partial L}{\partial u}h+\langle\nabla h, F\rangle \right)\, dxdy\\
&=\int_\Omega h\cdot \left(\frac{\partial L}{\partial u}-\mbox{div} F\right) \, dxdy +\int_\Omega\mbox{div}(h\cdot F)\, dxdy.
\end{split}
\end{equation*}
The Divergence Theorem allows us to rewrite the last integral as an integral over the boundary $\partial\Omega$. So, using  $h=0$ on $\partial\Omega$, we have
$$\int_\Omega\mbox{div}(h\cdot F)\, dxdy=\int_{\partial\Omega}h\cdot \langle F,{\bf n}\rangle=0,$$
where ${\bf n}$ is the unit outward-pointing normal of $\partial\Omega$.  As a consequence of the Fundamental Lemma of the calculus of variations, $u$ is an extremal if and only if  
$$\frac{\partial L}{\partial u}-\mbox{div} F=\frac{\partial L}{\partial u}-\left(\frac{\partial L}{\partial p}\right)_x-\left(\frac{\partial L}{\partial q}\right)_y=0.$$
By the   definition of $L$, this identity can be expressed as 
$$\sqrt{1+|\nabla u|^2}- \left((u+\lambda)\frac{ u_x}{\sqrt{1+|\nabla u|^2}}\right)_x- \left((u+\lambda)\frac{ u_y}{\sqrt{1+|\nabla u|^2}}\right)_y=0. $$
Rewriting this identity, we conclude that  an extremal $u$ of the variational problem satisfies  the Euler-Lagrange equation
$$ \left( \frac{ u_x}{\sqrt{1+|\nabla u|^2}}\right)_x+ \left( \frac{ u_y}{\sqrt{1+|\nabla u|^2}}\right)_y=\frac{1}{(u+\lambda)\sqrt{1+|\nabla u|^2}},$$
where $\lambda$ is a Lagrange multiplier. By translating the surface in the vertical position,  we can assume that $\lambda=0$, hence   
\begin{equation}\label{eq1}
\mbox{div}\frac{\nabla u}{\sqrt{1+|\nabla u|^2}}=\frac{1}{u\sqrt{1+|\nabla u|^2}}.
\end{equation}
Notice that we are only interested in smooth surfaces so, from \eqref{eq1}, we are only interested in functions $u(x,u)$ which are nowhere zero. It is also clear that if $u=u(x,y)$ satisfies \eqref{eq1}, then so does $-u(x,u)$. Thus, without loss of generality, we will only consider strictly positive functions $u(x,y)$ that satisfy \eqref{eq1}. In such a case, we  will  say that the surface $z=u(x,y)$ is a  {\it singular minimal surface}. This definition was first coined by Dierkes in \cite{di2} and motivated by  two reasons. First, because  the first integrand in  \eqref{jj}  degenerates if $u$ vanishes. Although as we said, our surfaces are smooth, in a more general context, one may admit singular solutions, that is,  $u=0$ somewhere. A simple example of a ``singular solution'' of   Equation \eqref{eqrot} is $u(r)=r$ in polar coordinates. This surface corresponds to a cone whose vertex is the origin and forming  a $45^{0}$  angle with respect to the rotational axis. Notice that this solution vanishes at $r=0$.

A second reason is that the left-hand side of \eqref{eq1} is a known quantity in differential geometry and  it coincides with the {\it mean curvature} $H$ of the surface $z=u(x,y)$. Minimal surfaces are those that satisfy $H=0$  everywhere.   In the one dimensional case (for functions of one variable), the solution of \eqref{eq1} is the catenary and for this, Equation \eqref{eq1}  is also known as the  {\it two-dimensional analogue of the catenary} (\cite{bht}). Indeed, Poisson already observed that if the function $u$ depends only on $x$, i.e.,   $u=u(x)$,   then \eqref{eq1} simplifies to 
$$\left(\frac{u'}{\sqrt{1+u'^2}}\right)'=\frac{1}{u\sqrt{1+u'^2}},\quad u'=\frac{du}{dx},$$
or equivalently,
\begin{equation}\label{cc}
\frac{u''}{1+u'^2}=\frac{1}{u},
\end{equation}
whose solution is the catenary $u(x)=a^{-1}\cosh(ax+b)$, $a\not=0$, $a,b\in\r$.  Of course,  the solutions of \eqref{eq1} are not minimal surfaces, but when we rotate the one dimensional solution (the catenary) with respect to the $x$-axis, we obtain the catenoid which  is a minimal surface. This explains that another reasonable name of a solution of \eqref{eq1} is   {\it symmetric minimal surface} (\cite{fs}).  A consequence of the derivation of the solutions of \eqref{eq1} in the one dimensional case is that the corresponding  cylindrical surface constructed with the catenary $u(x,y)=a^{-1}\cosh(ax+b)$ and repeated along the $y$-direction 
$$X(x,y)=(x,y,u(x))=(x,0,u(x))+y(0,1,0)$$
 is a singular minimal  surface whose rulings are all parallel to the $y$-axis. After inverting this surface, we again obtain the shape of  Gaud\'{\i}'s famous corridors.

We rewrite \eqref{eq1} in a form that will be used in the rest of the article.  Any surface of $\r^3$ is locally the graph of a function in one of the three coordinate planes of $\r^3$. We will assume that it is the graph over the $xy$-plane. Then locally $S=\{(x,y,u(x,y)):(x,y)\in\Omega\}$ for some smooth function $u$. We parametrize $S$ as 
$$X:\Omega\to\r^3,\quad X(x,y)=(x,y,u(x,y)).$$
We have pointed out that the left-hand side of \eqref{eq1} is just the mean curvature $H$ at any point $p\in S$. For the right-hand side, consider  the upward pointing unit normal vector field $N$ on $S$ which is orthogonal to the tangent plane. For the   parametrized surface $X(x,y)$,   the tangent plane is spanned by $\{X_x,X_y\}$, so $N$ can be computed by 
$$N=\frac{X_x\times X_y}{|X_x\times X_y|}=\frac{1}{\sqrt{1+|\nabla u|^2}}(-u_x,-u_y,1).$$
If $\vec{v}=(0,0,1)$ is the unit vector in the positive direction of the $z$-axis, then  $\langle X,\vec{v}\rangle=u$ and Equation  \eqref{eq1} can be expressed simply by 
\begin{equation}\label{hh}
H(p)=\frac{\langle N(p),\vec{v}\rangle}{\langle p,\vec{v}\rangle},\quad p\in S.
\end{equation}
Now the condition that the surface is a singular minimal surface is expressed free of coordinates. Many questions   are now open to us,  some of which have connections with the shape of a cupola. In this article we will deal with the following two aspects:
\begin{enumerate}
\item Given a closed curve $\Gamma$, does  the geometry of $\Gamma$ impose  restrictions to the shape of a singular minimal surface spanning $\Gamma$?
\item What is the shape of a surface of revolution that satisfies the singular minimal surface equation  \eqref{hh}?
\end{enumerate}

\begin{figure}[hbtp]
\begin{center}
\includegraphics[width=.32\textwidth]{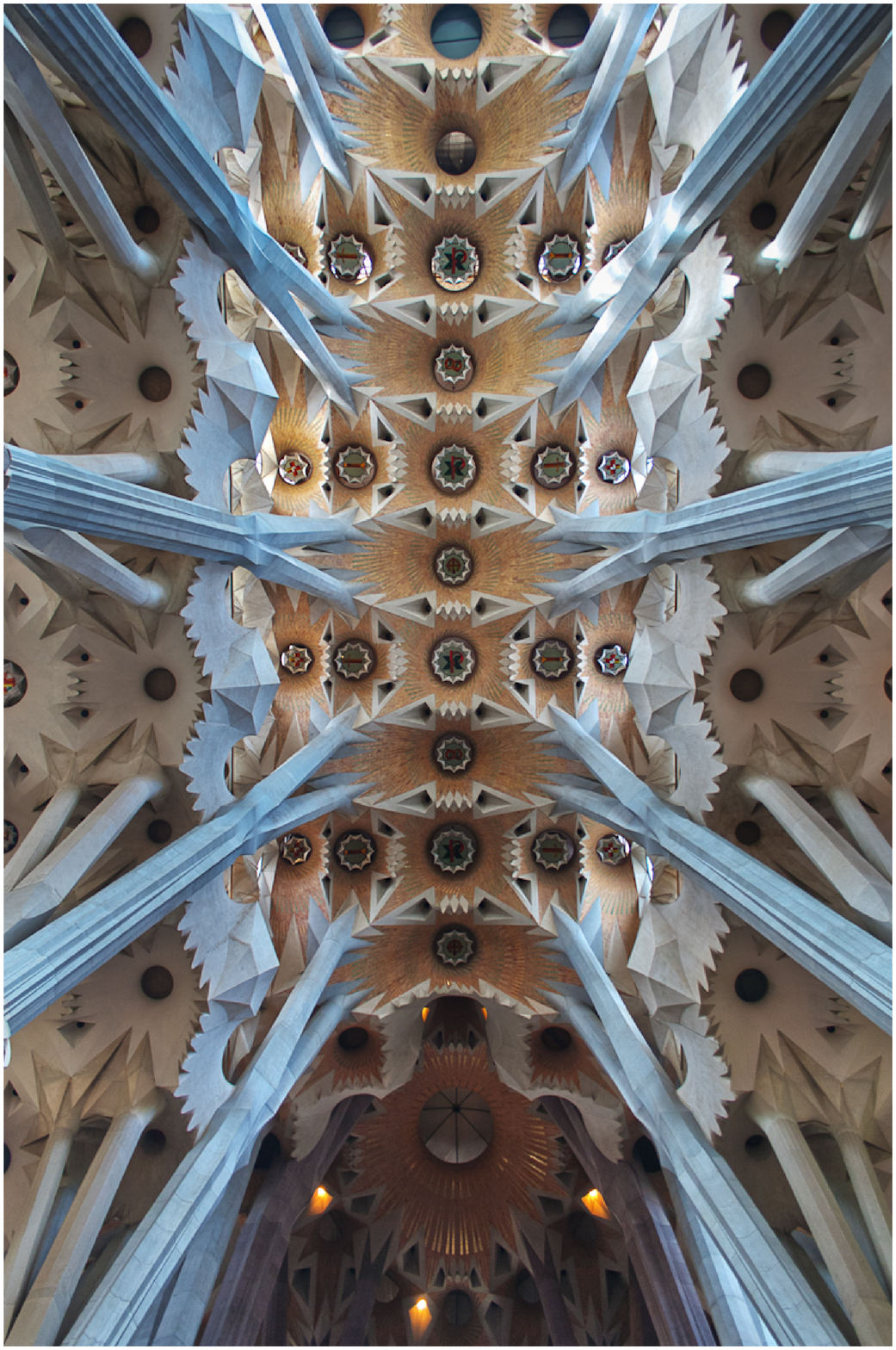}\qquad \includegraphics[width=.36\textwidth]{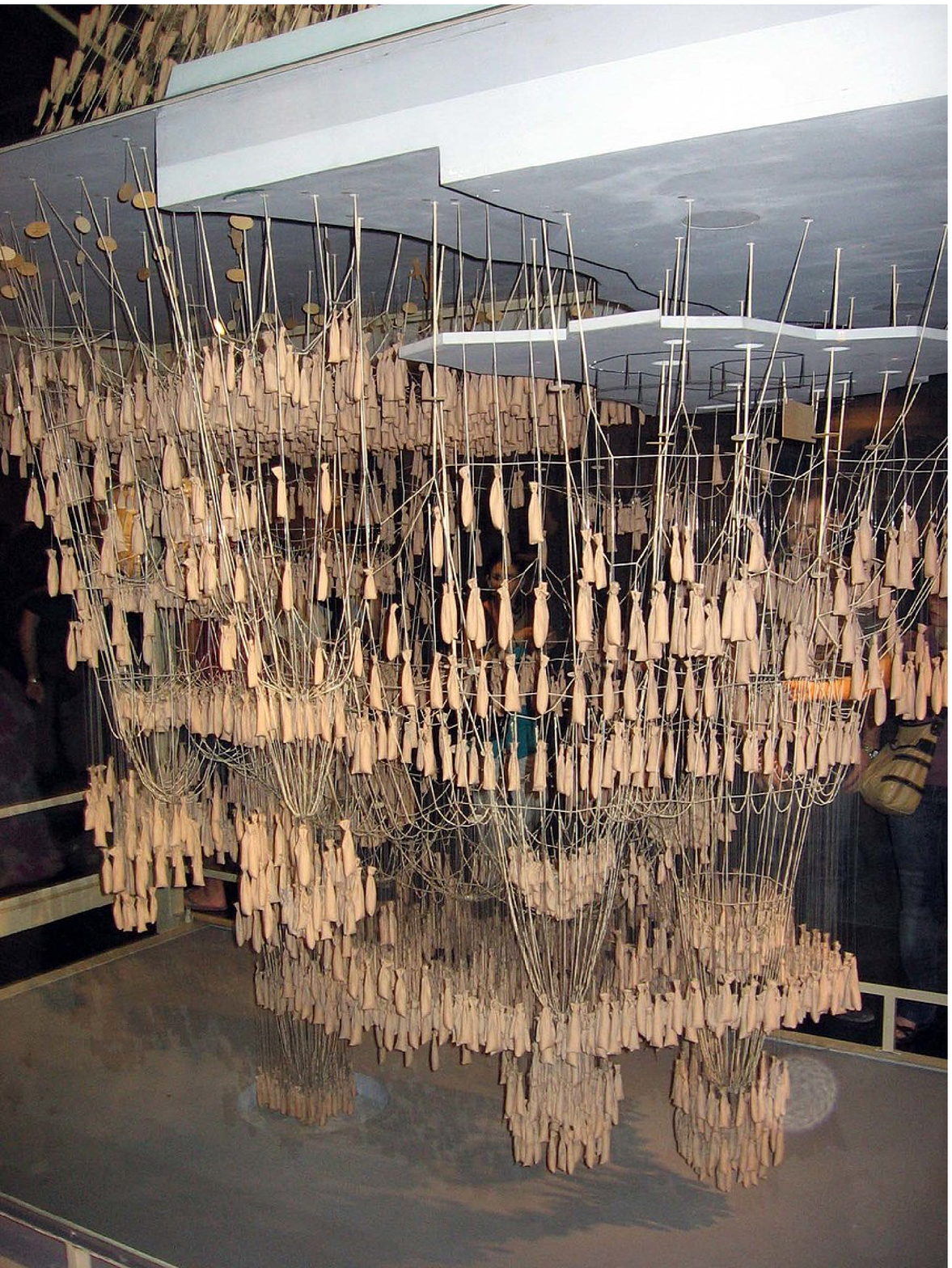}
\end{center}
\caption{Left: Dome on la Sagrada Familia, Barcelona, viewed from below. Right: model of the church of Colonia G\"{u}ell used by Gaud\'{\i} (Museum of Sagrada Familia). The first image   is licensed under the Creative Commons Attribution-Share Alike 3.0 Unported license at commons.wikimedia.org/wiki/File:Sagrada familia, Boveda principal.jpg. 
The second image is licensed under the GNU Free Documentation License at commons.wikimedia.org/wiki/File:Maqueta funicular.jpg.}\label{fig3}
\end{figure}

We return for the moment to historical considerations. As stated before,  equation \eqref{eq1} has been forgotten for some time. Independently, the shape of a cupola was addressed by Antonio Gaud\'{\i}. He was particularly interested in the shape of a suspended surface.   For the construction of his unfinished work on the basilica known as the Sagrada Familia (Sacred Family), see Figure \ref{fig3}, left, he wanted to reproduce the shape  of these surfaces. This was important to him because of his own style of reproducing forms from nature. So,  Gaud\'{\i}   designed the structure of  a dome  by suspending loads from wires that simulated the different arches and pillars upside down, as can be seen in Figure \ref{fig3}, right. Many years later, Frei Otto again reproduced  this design in the Institute for Lightweight Structures at the University of Stuttgart (\cite{ot2}).

We present some properties of the singular minimal surfaces.
\begin{enumerate}

\item {\it The set of singular minimal surfaces is invariant by rigid motions that fix the vertical direction and whose translation vectors are horizontal. They are also invariant by dilations from any point of the $xy$-plane.} Indeed,  let  $M\colon\r^3\to\r^3$ be a rigid motion, with $M(p)=A  p+\vec{b}$, where $A\in O(3)$ is a linear isometry and $\vec{b}\in\r^3$ is the translation vector. Let $\tilde{S}=M(S)$ and denote $\tilde{p}=M(p)=Ap+\vec{b}$ for $p\in S$.  At corresponding points $p$ and $\tilde{p}$,  the mean curvature coincide and $N(\tilde{p})=AN(p)$. Thus  \eqref{hh} becomes 
$$H(\tilde{p})=\frac{\langle A^{-1}N(\tilde{p}),\vec{v}\rangle}{\langle A^{-1}\tilde{p},\vec{v}\rangle-\langle A^{-1}\vec{b},\vec{v}\rangle}=\frac{\langle N(\tilde{p}),A\vec{v}\rangle}{\langle \tilde{p},A\vec{v}\rangle-\langle A^{-1}\vec{b},\vec{v}\rangle}.$$
Since we want to keep the vertical direction for  gravity, we need that $A\vec{v}=\vec{v}$, concluding $\langle\vec{b},\vec{v}\rangle=0$, too. Examples of these motions are rotations about a vertical straight line, symmetries about vertical planes, and horizontal translations (recall that the vertical translations were already used in the elimination of $\lambda$ in the derivation of \eqref{eq1}). The proof for dilations is straightforward because we can again assume that, after a horizontal translation, the dilation is expressed by $h(p)=\lambda p$, $\lambda>0$. In such a case, if $\tilde{p}=h(p)$, then $H(\tilde{p})=H(p)/\lambda$ and $N(\tilde{p})=N(p)$.

\item Writing the surface as $z=u(x,y)$, we see that {\it the function $u$ has no local maximum at the interior points of $\Omega$}. Indeed, if $q=(x_0,y_0)\in\Omega$ is a local maximum, then $\nabla u(q)=(0,0)$ and $u_{xx}(q)\leq 0$ and $u_{yy}(q)\leq 0$. If we expand out \eqref{eq1}, we have 
\begin{equation}\label{eqpo}
(1+u_y^2)u_{xx}-2u_x u_y u_{xy}+(1+u_x^2)u_{yy}=\frac{1}{u}(1+u_x^2+u_u^2).
\end{equation}
At $q$, this identity reduces to
$$(u_{xx}+u_{yy})(q)=\frac{1}{u(q)}>0,$$
which is not possible. Notice that \eqref{eqpo} coincides with the equation of Figure \ref{figpo}, up to the constants $g$ and $\epsilon$.

\item As a consequence,  {\it if the boundary $\Gamma$ of $S$ is contained in a horizontal plane and $S$ is compact, then $S$ lies below that plane.}
\item Suppose the function  $\varphi\colon \partial\Omega\to\r$ defining the boundary curve $\Gamma$ is smooth. We prove that {\it if $u$ satisfies \eqref{eq1}, then}
\begin{equation}\label{al}
\frac{\mbox{area}(\Omega)}{\mbox{length}(\partial\Omega)}<\max_{\partial\Omega} \varphi.
\end{equation}
This gives a necessary condition in terms of the geometry of the boundary curve $\Gamma$ for the existence of a singular minimal surface spanning $\Gamma$. The proof of \eqref{al} is as follows. A simple computation yields 
\begin{equation*}
\begin{split}
\mbox{div}\left(u\cdot \frac{\nabla u}{\sqrt{1+|\nabla u|^2}}\right)&=\frac{|\nabla u|^2}{\sqrt{1+|\nabla u|^2}}+u\cdot \mbox{div}\left( \frac{\nabla u}{\sqrt{1+|\nabla u|^2}}\right)\\
&=\sqrt{1+|\nabla u|^2}.
\end{split}
\end{equation*}
Integrating over $\Omega$ and using the Divergence Theorem, we obtain
$$\int_\Omega\sqrt{1+|\nabla u|^2}\ dx dy= \int_{\partial\Omega}\varphi\frac{  \langle {\bf n},\nabla u\rangle }{\sqrt{1+|\nabla u|^2}}.$$
  The left-hand side in the above identity is the area of $S$. On the other hand, since $|\langle {\bf n},\nabla u\rangle|\leq|{\bf n}| |\nabla u|=|\nabla u|$, we deduce 
$$\mbox{area}(S)< \int_{\partial\Omega} \max_{\partial \Omega}\varphi=(\max_{\partial \Omega} \varphi)\cdot \mbox{length}(\partial\Omega).$$
Hence, because  $\mbox{area}(\Omega)<\mbox{area}(S)$,  the inequality  \eqref{al} holds, as claimed. Here we point out that Nitsche already gave an upper bound   $\mbox{area}(S)< A(\Gamma)$, $A(\Gamma)$ depending only on $\Gamma$, when $S$ is a rotational singular minimal surface and $\Gamma$ is a horizontal circle (\cite{ni}).

\item As a consequence of \eqref{al}, the existence of a solution to the Dirichlet problem associated to \eqref{eq1} with boundary conditions $u=\varphi$ on $\partial\Omega$ is not assured for general $\varphi$. On the other hand, it is also not known under what conditions one has uniqueness of solutions for the Dirichlet problem, or whether a solution is in fact a minimizer of the variational problem. See \cite{lo3,lo4}.
\end{enumerate}

We conclude this section with an expected property that requires difficult techniques beyond the scope of this article.   Suppose that the boundary $\Gamma$ is just a circle contained in a horizontal plane $\Pi$. If $S$ is a compact singular minimal surface spanning $\Gamma$, does $S$ inherit the axisymmetric shape of $\Gamma$? The answer is yes if we assume that $S$ is a surface without self-intersections (as is the case for graphs). Firstly, by property (3), $S$ lies below $\Pi$. Now an argument  due to Alexandrov (\cite{al}) using reflection across vertical planes together with the Maximum Principle, proves that given any vertical plane $P$, there is another parallel plane to $P$ such that $S$ is invariant by reflections across that plane.  Doing this for every vertical plane, one concludes that, indeed, the surface is rotationally symmetric about a vertical line (\cite{lo2}).

\begin{theorem}\label{t2} Let $S$ be a compact singular minimal surface without self-intersections. If the boundary is a horizontal circle, then $S$ is a surface of revolution about a straight line parallel to the $z$-axis.
\end{theorem}

 \section{Rotational cupolas.}\label{sec3}
 
\begin{figure}[hbtp]
\begin{center}
\includegraphics[width=.4\textwidth]{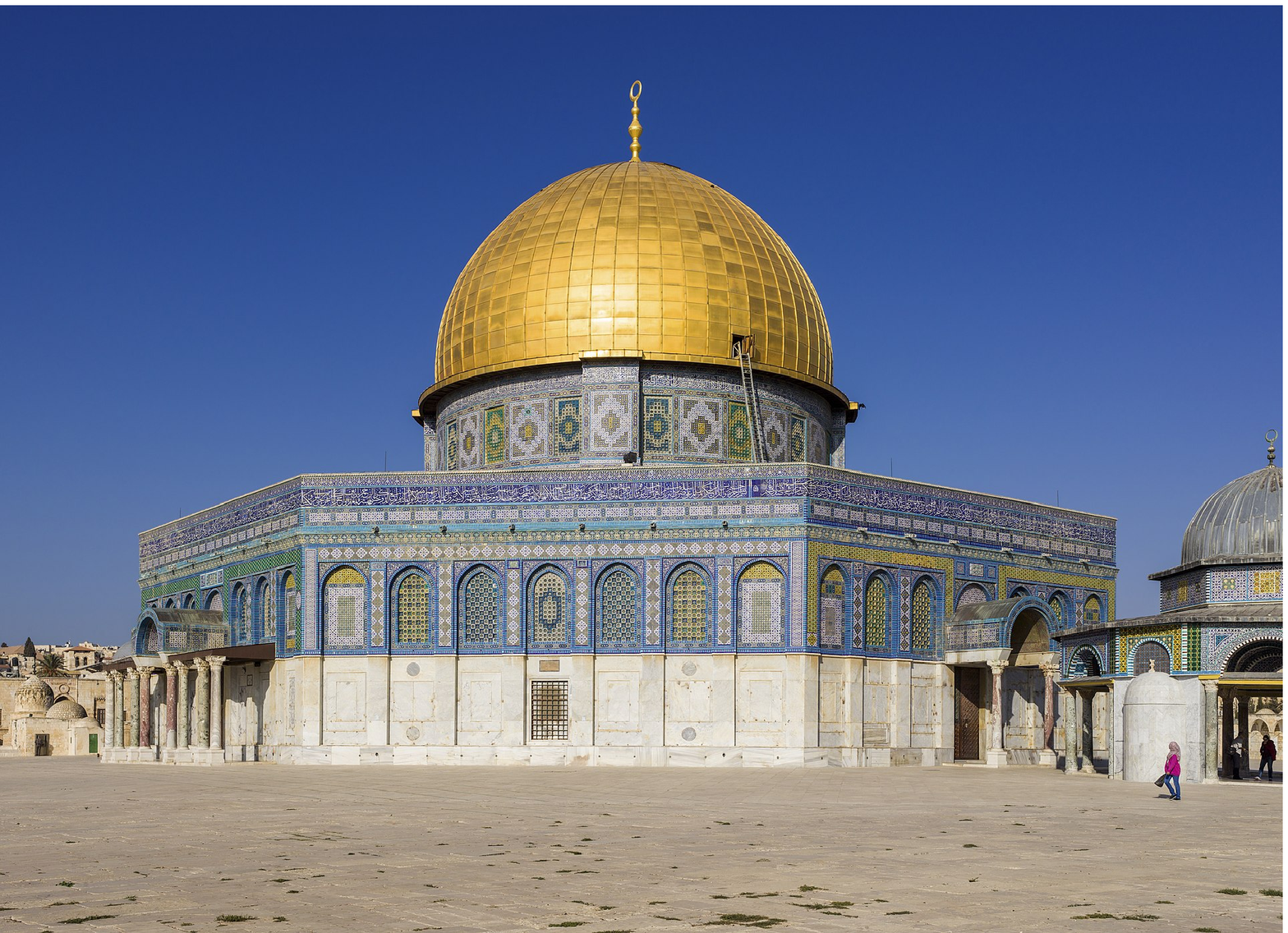}\qquad \includegraphics[width=.5\textwidth]{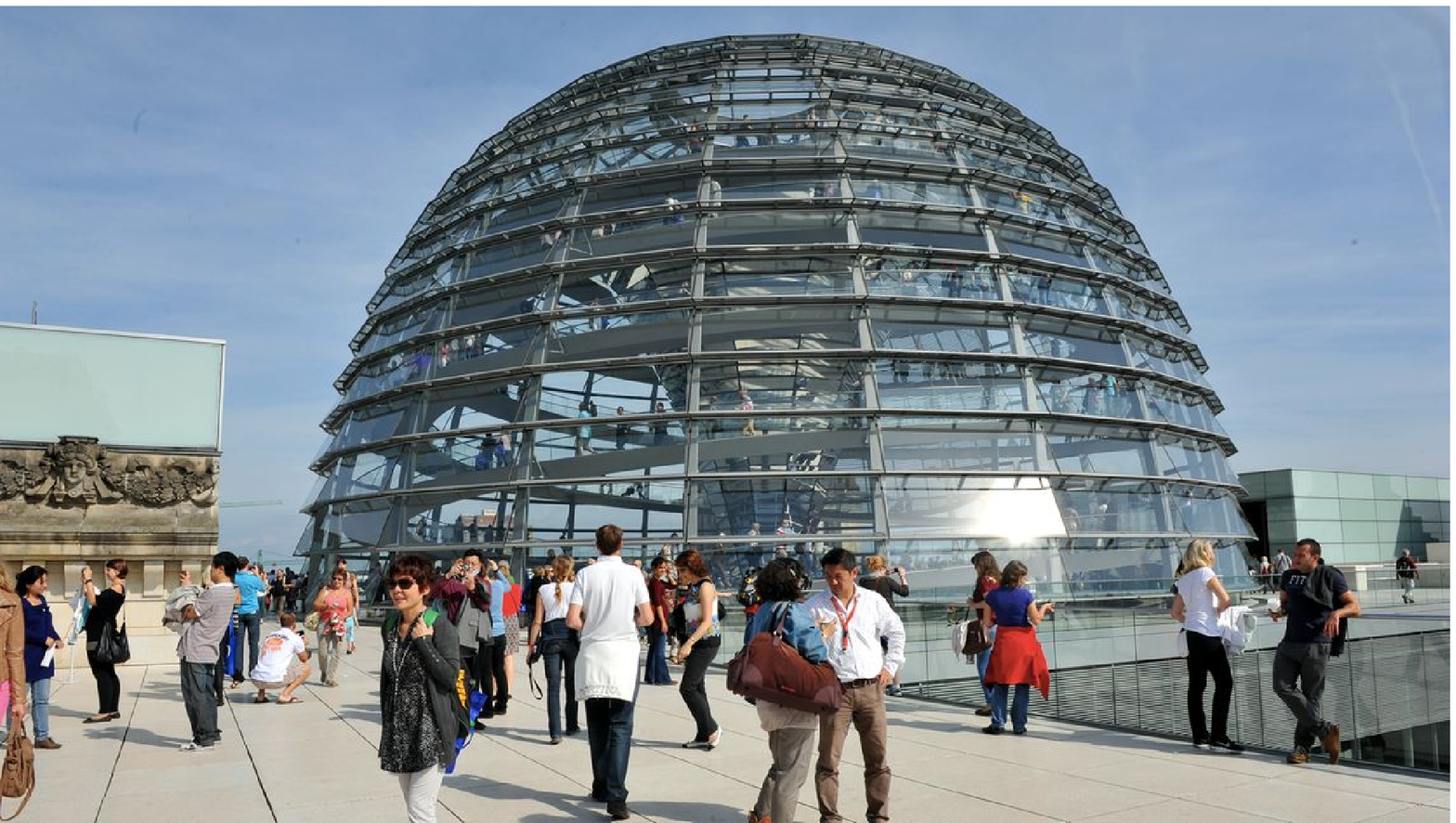}
\end{center}
\caption{Cupolas that look like surfaces of revolution: the Dome of the Rock at Jerusalem (left) and the Reichstag   in Berlin by Norman Foster (right). The first image is licensed under the Creative Commons Attribution-Share Alike 4.0 Unported license at commons.wikimedia.org/wiki/File:Jerusalem-2013(2)-Temple Mount-Dome of the Rock (SE exposure).jpg (Andrew Shiva). The source of the second image is  \url{https://www.bundestag.de/en/visittheBundestag/dome/registration-245686},  $\textcopyright$ Deutscher Bundestag/Neuhauser.}\label{fig4}
\end{figure}

 Let us come back to our initial problem of the construction of cupolas. The common idea to build a cupola is that its shape is modeled by a surface of revolution whose rotation axis is a vertical line. Even in this case, the real construction of a rotational cupola never occurs because architects use `discrete' methods of construction. So, the base of the cupola is never a circle, but it is a `discrete' circle formed by  a union of rectilinear segments adopting circular shape. In fact, the cupolas of Figure \ref{fig1} are not surfaces of revolution: their shape is invariant by a finite group of rigid motions, which coincides with the number of arches connecting the top of the cupola to its base. In the case of the cupola of Brunelleschi, this number is $8$, while it is 16 in that of Michelangelo.   Other cupolas whose shapes better resemble surfaces of revolution are shown in Figure \ref{fig4}.

Although we know that cupolas are surfaces of revolution, from a mathematical perspective it is not clear that their rotation axes must be parallel to the direction of gravity. We investigate this question.   To facilitate the computations, we suppose that the rotation axis of the surface $S$  is the $z$-axis but the direction of gravity is indicated  by the  direction $\vec{v}=(a_1,a_2,a_3)$ with $|\vec{v}|^2=1$.  All points $p=(x,y,u(x,y))\in S$ that lie at the same horizontal plane, are circles centered at the $z$-axis of radius $r=\sqrt{x^2+y^2}$. Thus $u$ is a radial function $u=u(r)$. Let us parametrize $S$ by introducing   polar coordinates $x=r\cos\theta$, $y=r\sin\theta$,  
\begin{equation}\label{ss}
X(r,\theta)=(r\cos\theta, r\sin\theta, u(r)).
\end{equation}
We express \eqref{hh} in terms of the derivatives of $u$ with respect to $r$. A change of variables transforms the left-hand side of \eqref{eq1} (equivalent to  the mean curvature $H$ in \eqref{hh}) into
\begin{equation}\label{uu}
\frac{u'(1+u'^2)+r u''}{r(1+u'^2)^{3/2}}.
\end{equation}
We now compute the right-hand side in \eqref{hh}. The unit normal vector field of $S$  is 
$$N=\frac{X_r\times X_\theta}{|X_r\times X_\theta|}= \frac{1}{\sqrt{1+u'^2}}(-u'\cos\theta,-u\sin\theta,1).$$
Since $\langle X,\vec{v}\rangle=a_1 r\cos\theta+a_2r\sin\theta+a_3 u$,   Equation \eqref{hh} is 
$$\frac{u'}{r}+\frac{u''}{1+u'^2}=\frac{-a_1u'\cos\theta-a_2u'\sin\theta+a_3 }{a_1 r\cos\theta+a_2r\sin\theta+a_3 u}.$$
After some manipulations, this equation can be written as $A(r)\cos\theta+B(r)\sin\theta+C(r)=0$, where 
\begin{equation*}
\begin{split}
A(r)&=a_1r\left( u'(1+u'^2)+ru''+   u'(1+u'^2)\right)\\
B(r)&=a_2r\left(u'(1+u'^2)+ru''+ u'(1+u'^2)\right)\\
C(r)&=a_3\left(u(u'(1+u'^2)+ru'')-r (1+u'^2)\right).
\end{split}
\end{equation*}
Since the functions $\{1,\cos\theta,\sin\theta\}$ are linearly independent, the functions $A$, $B$ and $C$ must vanish in their domain. One case is that $a_1=a_2=0$. Then the direction of gravity is parallel to the rotation axis and, in addition, $C=0$ becomes  
 \begin{equation}\label{eqrot}
\frac{u''}{1+u'^2}=\frac{1}{u}-\frac{u'}{r}.
\end{equation}
 Suppose now that $a_1\not=0$ (resp. $a_2\not=0$). Then we deduce from $A=0$ (resp. $B=0$), 
\begin{equation}\label{eq2}
\frac{u''}{1+u'^2}=-\frac{2u'}{r}.
\end{equation}
 For the equation $C(r)=0$, we distinguish two subcases. If  $a_3=0$, then $\vec{v}$ is orthogonal to the axis of rotation. If $a_3\not=0$, then   $u(u'(1+u'^2)+ru'')-r (1+u'^2)=0$ and combining with \eqref{eq2}, we deduce $uu'=-r$. Thus 
 $u(r)=\sqrt{r^2+c}$, for some constant $c$. However, this function does not satisfy \eqref{eq2}. This establishes the following theorem which is now written when the direction of the gravity is given, as usual, by the vertical axis (\cite{lo1}). 
 
 \begin{theorem} \label{t1}
 If a surface of revolution   is a singular minimal surface, then the axis of rotation is vertical  or the axis of rotation is contained in the plane $z=0$. 
 \end{theorem}

The second case is striking because we have discovered a model of a rotational cupola whose rotation axis is horizontal! We separate the two cases and, in this section,  we focus on the case where the rotation axis is parallel to the force of gravity. Here $\vec{v}$ in the proof of the theorem is actually the vertical direction of $\r^3$ and $u$ satisfies \eqref{eqrot}. From the standard theory of ordinary differential equations, the solution of the ordinary differential equation \eqref{eqrot} is obtained once we give initial conditions
\begin{equation}\label{eqrot2}
u(r_0)=u_0,\quad u'(r_0)=\bar{u}_0,\quad r_0>0, u_0>0.
\end{equation}
Let us observe  that  \eqref{eqrot} is singular at $r=0$, so $r_0$ must be positive. However, keeping in mind the shape of cupolas, our interest is that a solution meets the rotation axis. So we want to know if the solution $u$ can be prolongated until $r=0$. This question is problematic. It is possible that under some initial conditions in \eqref{eqrot2}, the solution does not meet the $z$-axis (see the example in Remark \ref{re1} below).  In such a case, after rotating the graphic of $u=u(r)$ about the $z$-axis, we would obtain  a cupola with a ``hole'' at the top.   However we are interested in those solutions whose initial conditions  in \eqref{eqrot2} occur at $r_0=0$. An argument using the Banach Fixed Point Theorem proves the existence of a solution with $u(0)=u_0>0$ (\cite{lo1}).  
 In such a case, the intersection of the surface with the rotational axis must be orthogonal by smoothness of the surface. This is equivalent to $u'(0)=0$.  

Rotational singular minimal surfaces whose axis is vertical have been studied in the literature: see, for example, \cite{di1,di2}. Recall that in Section \ref{sec2} we showed the singular solution $u(r)=r$, a cone with vertex at the origin. In this case, after inverting the surface, the shape of the cupola looks like a Native American teepee. 

\begin{figure}[hbtp]
\begin{center}
\includegraphics[width=.4\textwidth]{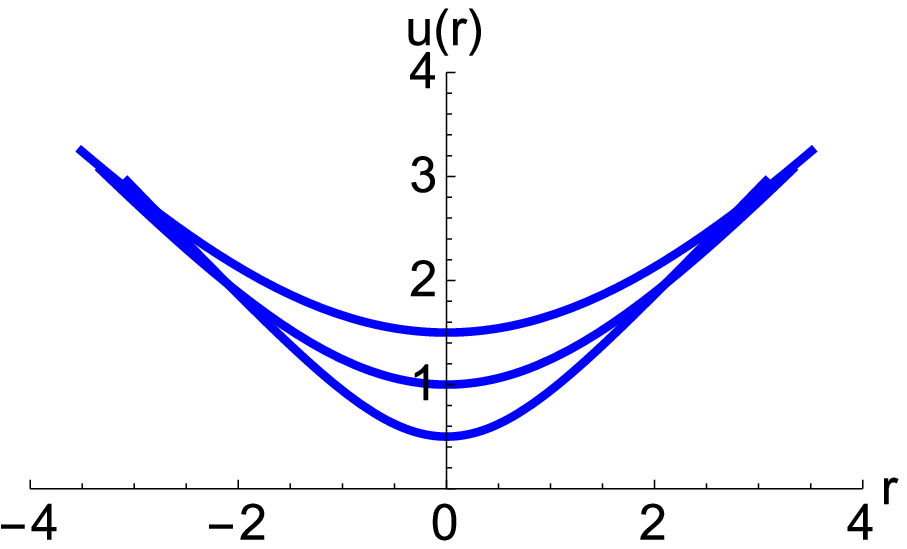}\qquad \includegraphics[width=.4\textwidth]{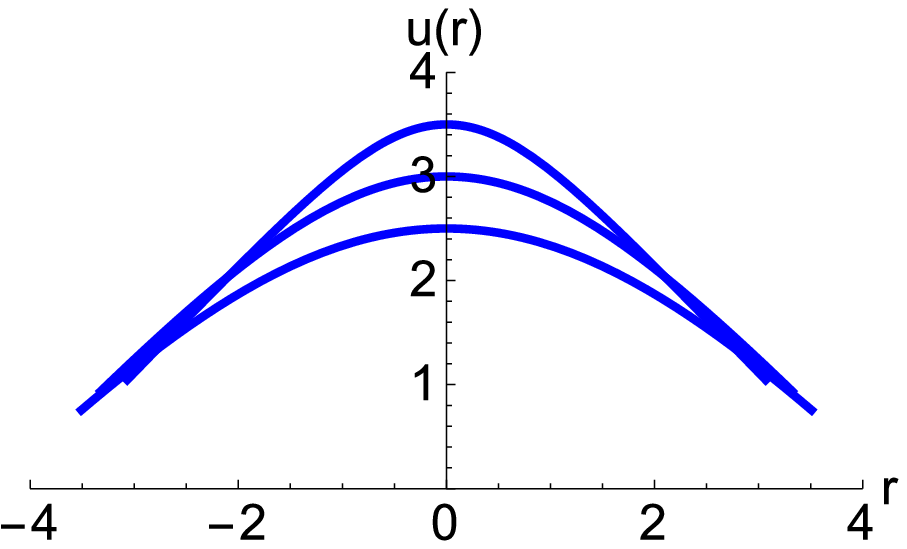}
\end{center}
\caption{Left: numerical solutions of \eqref{eqrot}-\eqref{eqrot2}, where $u'(0)=0$ for different values of $u(0)$:  $0.5$, $1$ and $1.5$. Right: the solutions viewed as cupolas after a symmetry about the horizontal line $z=2$.}\label{fig5}
\end{figure}

In Figure \ref{fig5}, left, we show, using {\sc Mathematica},  some numerical solutions of \eqref{eqrot}-\eqref{eqrot2} when $r_0=0$, $u'(0)=0$, and  for different values of $u_0$.  All these curves will give shapes of domes once we invert them as  Figure \ref{fig5}, right, shows.

\begin{remark} \label{re1} 
If $r_0>0$ in \eqref{eqrot2}, then the standard theory of ODE's ensures the existence and uniqueness of solutions. In such a case, the maximal domain of the solution around $r_0$ may not reach the value $0$, that is, the solution may not meet the rotation axis. This happens  when we choose $u'(r_0)=0$ for a fix value $r_0>0$.  Indeed, if the domain of $u$ contains the value $0$, we know that $u'(0)=0$. From \eqref{eqrot}, $u''(r_0)=1/u(r_0)>0$, so $r=r_0$ is a strict local minimum. We now see that $r=0$ is another strict local minimum. By   L'H\^{o}pital's rule, letting $r\to 0$ in \eqref{eqrot}, we get
$$u''(0)=\frac{1}{u(0)}-\lim_{r\to 0}\frac{u'(r)}{r}=\frac{1}{u(0)}-\lim_{r\to 0}\frac{u''(r)}{1}=\frac{1}{u(0)}-u''(0).$$
Hence $u''(0)=1/(2u(0))>0$. Thus the function $u$ restricted to the interval $[0,r_0]$ must have a local maximum at some point $c\in (0,r_0)$, which must also be a local maximum of $u=u(x,y)$. This contradicts property (2) of Section \ref{sec2}.  In Figure \ref{fig52} we show an example of a rotational surface that does not meet the $z$-axis.
\begin{figure}
\begin{center}
\includegraphics[width=.3\textwidth]{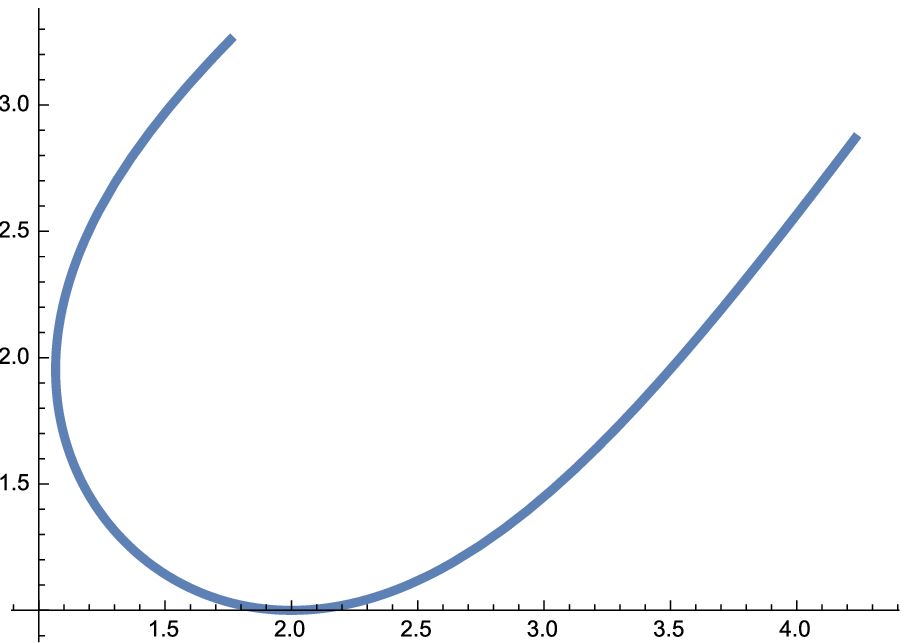}\qquad \includegraphics[width=.4\textwidth]{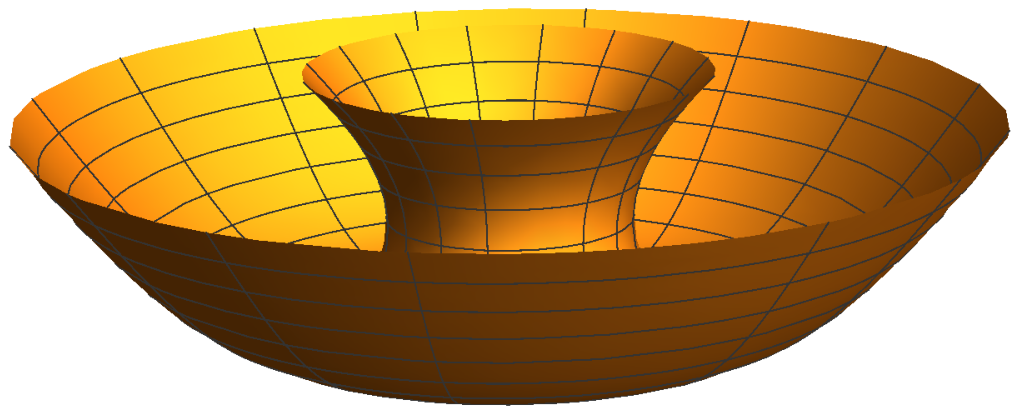}
\end{center}
\caption{Left: a solution of  \eqref{eqrot}-\eqref{eqrot2}, where $u(2)=1$ and $u'(2)=0$. Right:  the corresponding rotational surface.}\label{fig52}
\end{figure}
\end{remark}

\section{A new design for a roof.}\label{sec4}

In this section we present a new design for a cupola using a   surface of revolution but, contrary to common sense, the rotational axis will be  horizontal. In this case, we feel it is better to refer to the surface as a `roof' rather than a cupola. Thus, we turn our attention to the singular minimal surfaces given in Theorem \ref{t1} whose rotation axis is included in the plane $z=0$. First, we need to change the coordinates in the proof of Theorem \ref{t1} because  here we assumed that the rotation axis was the $z$-axis and the direction $\vec{v}$ was $(a_1,0,0)$ or $(0,a_2,0)$.  Without loss of generality, we suppose $\vec{v}=(1,0,0)$ and consider the positively oriented rigid motion $M\colon\r^3\to\r^3$  determined by the transformation 
$$M\colon \left\{\begin{array}{l}(1,0,0)\mapsto (0,0,1)\\ (0,1,0)\mapsto (0,-1,0)\\ (0,0,1)\mapsto (1,0,0).\end{array}\right.$$
The  surface of revolution  in \eqref{ss}  changes to $M\circ X(r,\theta)=(u,-r\sin\theta,r\cos\theta)$. On the other hand, we know that $u$ satisfies \eqref{eq2}.  We make a new change of variables interchanging the roles of  $u$ and $r$. Then  the parametrization of the surface is 
\begin{equation}\label{x5}
X(r,\theta)=(r,-u(r)\sin\theta,u(r)\cos\theta),
\end{equation} 
and   \eqref{eq2} is now
\begin{equation}\label{eq3}
\frac{u''}{1+u'^2}=\frac{2}{u}.
\end{equation}
Notice that this equation looks like the equation \eqref{cc} of the catenary with the only difference being that now the numerator in the right-hand side is $2$. For this reason, we call  the solutions of \eqref{eq3} {\it $2$-catenaries}.  For non-constant solutions, we multiply  by $u'$ and integrating, we find 
\begin{equation}\label{fin}
u'=\pm \sqrt{c^2 u^{4}-1}, \quad c\not=0.
\end{equation}
 In particular, differentiating with respect to $r$, and using \eqref{fin}  
\begin{equation}\label{ac3}
u''(r)=  \pm 2c^2\frac{u^3u'}{\sqrt{c^2u^4-1}}=2c^2u(r)^3.
\end{equation}
 The differential equation \eqref{fin} is known in the literature as an Emden-Fowler type equation (\cite{pz}). The generating curve $u=u(x)$ is contained in the $xz$-plane after we replace the  variable $r$ with $x$. The properties of $u$ are the following:
 \begin{enumerate}
 \item The function $u$  has only one critical point. Without loss of generality, we can assume that this point is $x=0$. At $x=0$, $u$ has a global minimum. The value of $u$ in $x=0$ is $1/\sqrt{c}$ by \eqref{fin}. 
  \item The function $u$ is symmetric about the $z$-line. The maximal domain  of $u$ is a bounded interval $(-a,a)$ and $\lim_{x\to\pm a}u(x)=\infty$. 
\item The function $u$ is convex thanks to \eqref{ac3}.
 \end{enumerate}
 
 If we were to build the roof rotating the curve $u=u(x)$ around the $x$-axis, the projection of the roof would be included in the horizontal strip $\Omega=\{(x,y,0)\in\r^3: -a\leq x\leq a\}$ and its walls, or its skeleton structure, would be near vertical at far away points. 
 
\begin{figure}[hbtp]
\begin{center}
\includegraphics[width=.3\textwidth]{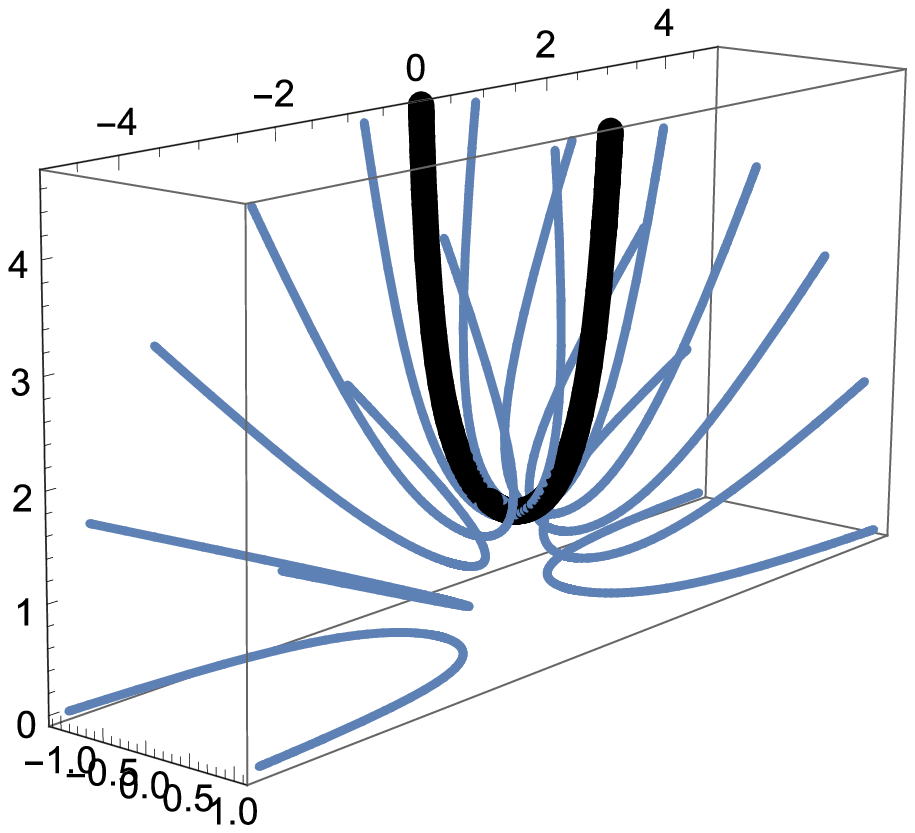}\quad \includegraphics[width=.32\textwidth]{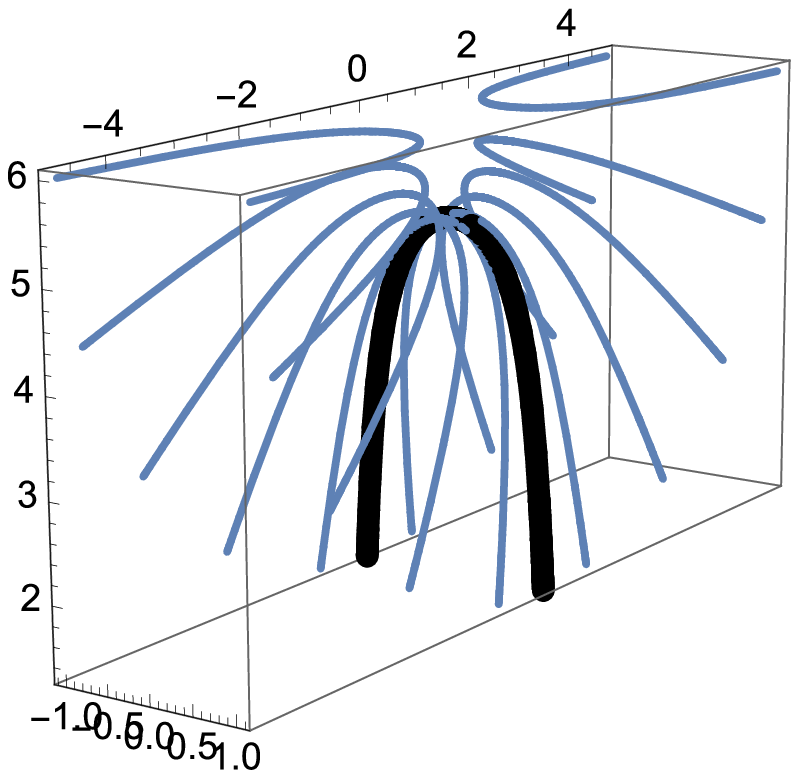}\quad \includegraphics[width=.28\textwidth]{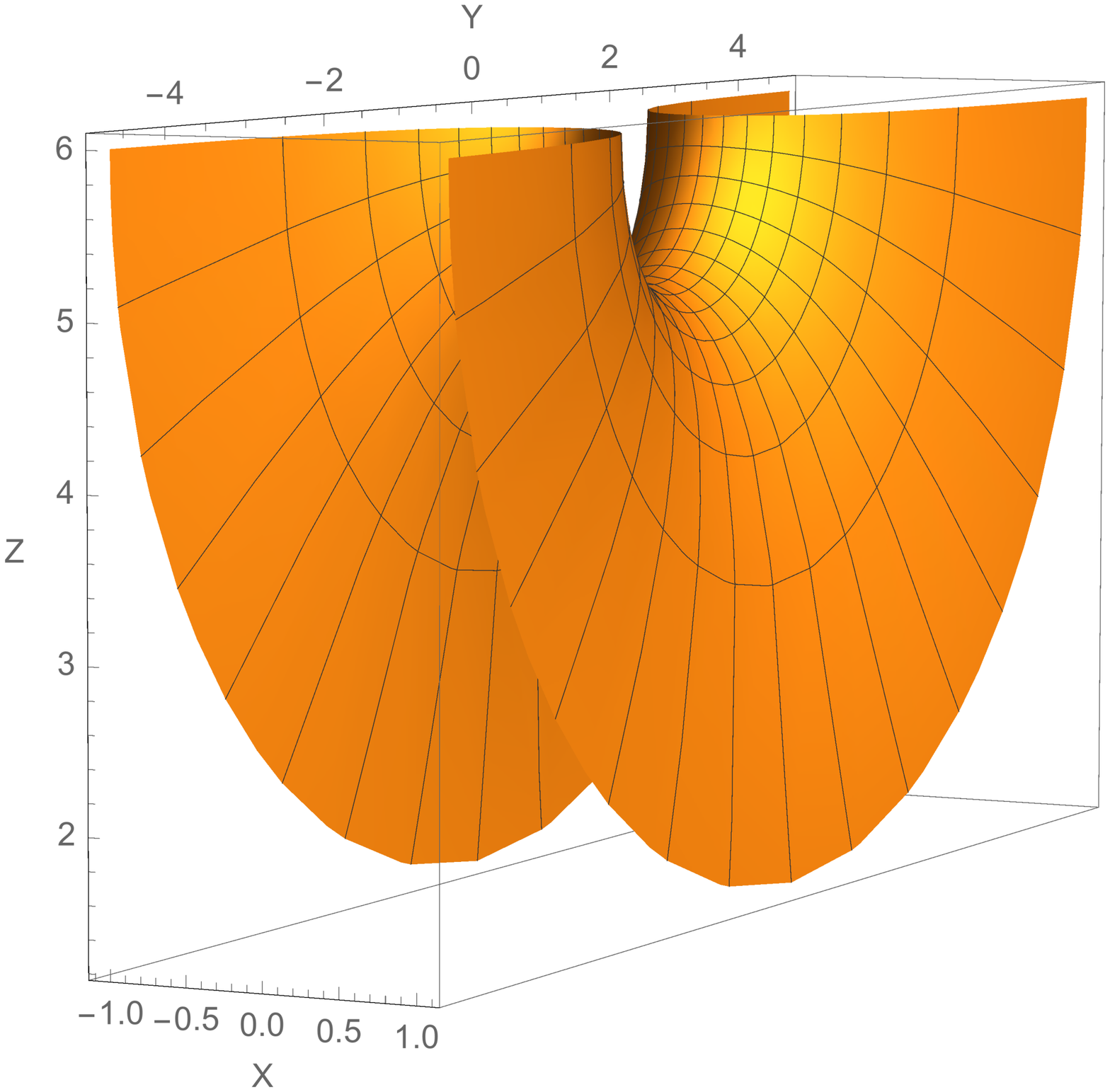}
\end{center}
\caption{Left: the $2$-catenary (bold) and its rotations about the $x$-axis. Middle: the same curves after inverting. Right: the rotational cupola. }\label{fig6}
\end{figure}

 We  plot numerical solutions of \eqref{eq3}   using {\sc Mathematica}. For this, we consider initial conditions 
$$ u(0)=1,\quad \ u'(0)=0.$$
 The maximal domain    $(-a,a)$ occurs for the value $a\approx  1.31102$. 
 The surface  is tangent to the vertical planes of equations $x=-a$ and $x=a$. When we rotate about the $x$-axis, the lower half of the surface is located in the halfspace $z<0$ which cannot be considered. In Figure \ref{fig6}, left, we plot the generating curve (thick) and the corresponding rotations of this curve for angles  $\theta\in (-\pi/2,\pi/2)$ (thin). In Figure \ref{fig6}, middle, we invert with respect to the horizontal plane having equation $z=3$ and on the right, we show the roof modeled by the surface. Both vertical walls are supporting the roof. Notice that all points of the surface are saddle points.

\begin{figure}[hbtp]
\begin{center}
\includegraphics[width=.32\textwidth]{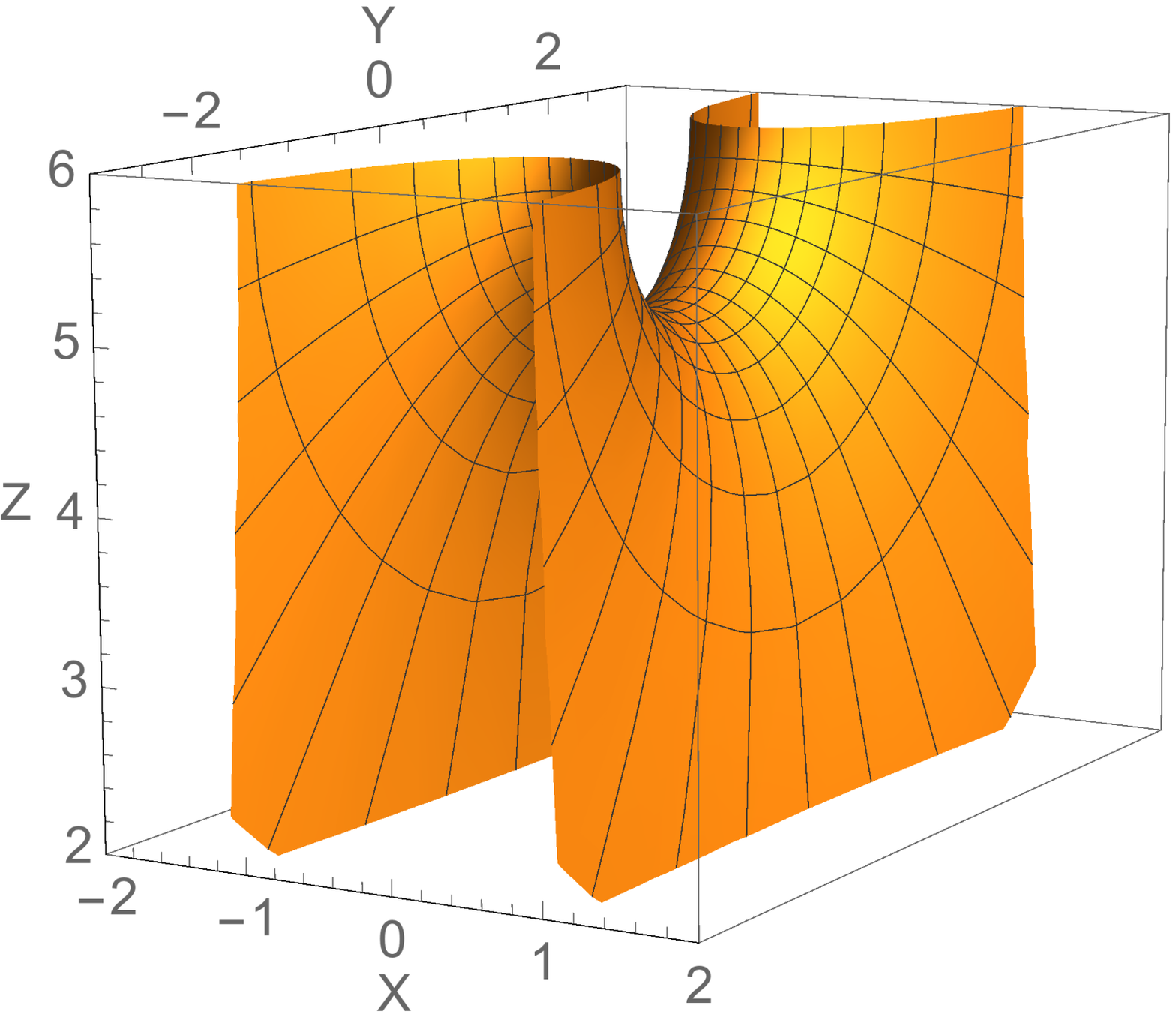}\quad \includegraphics[width=.3\textwidth]{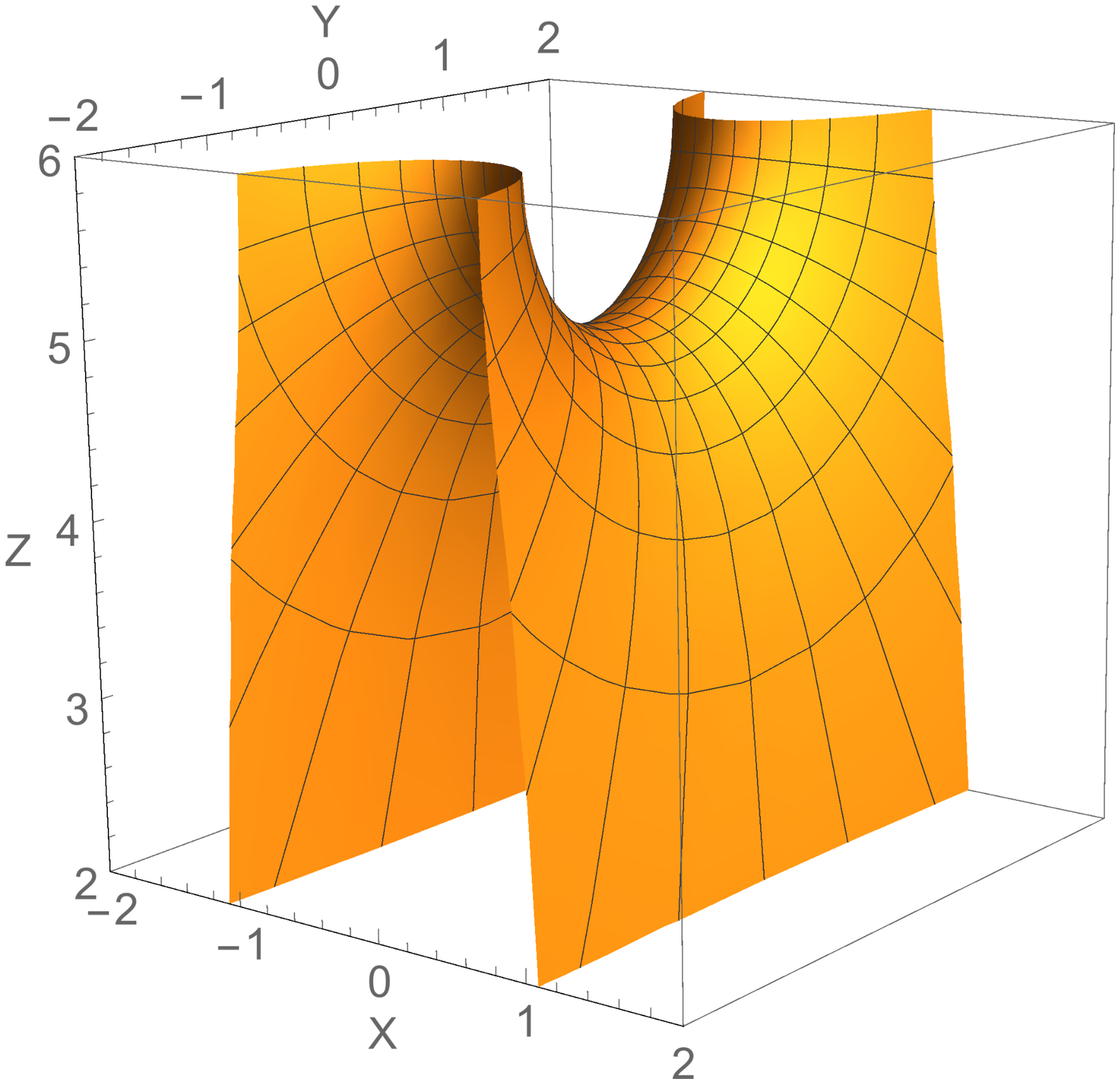}\quad \includegraphics[width=.27\textwidth]{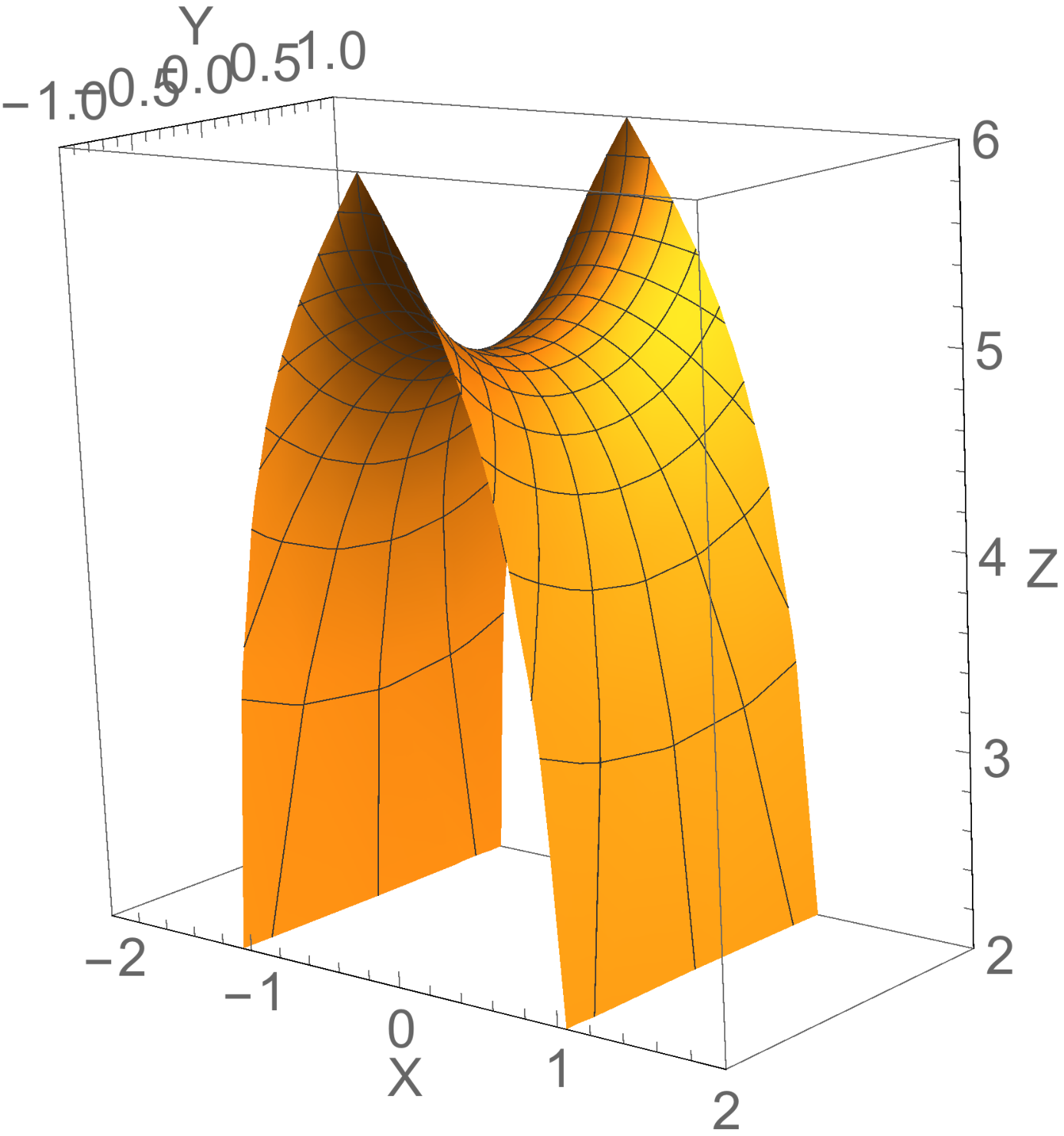}
\end{center}
\caption{Rotational cupolas by cutting along vertical planes parallel to the $xz$-plane. The values for $y$ are varying in the interval $[-3,3]$ (left), $[-2,2]$ (middle) and $[-1,1]$ (right).  }\label{fig7}
\end{figure}

In the last surface in Figure \ref{fig6}, the roof  does not cover the entire corridor (the strip $\Omega$). We have reduced the size of the roof along the $y$-axis as  Figure \ref{fig7} shows.

\section{Outlook and Conclusions.} Motivated by the shape of a catenary, we have deduced the differential equation governing surfaces suspended by their own weight and discussed some of their properties. Singular minimal surfaces can be models for cupolas, at least under the simple hypothesis of compression. Hence, light structures can be constructed in architecture imitating the shape of these surfaces.

In reality, these surfaces may be difficult to produce on  a large scale. However, it is remarkable that there are singular minimal surfaces that are surfaces of revolution about a horizontal axis. Thanks to these surfaces, we have presented   a novel structural shape of a roof in Figures \ref{fig6} and \ref{fig7} which can be regarded in the context of the so-called ``funicular shape'' in architecture. The two families of parametric curves in this surface show the visual design of a skeleton that opens up towards its border, increasing its beauty.  And, as has been justified, its shape is `natural' in the sense that loads and tensions act tangentially on the roof, giving solidity and stability to the construction.

Gaud\'{\i} used principles from the natural sciences in his architecture, generating interest in the design of structures by observing the effect of weight.  The idea to create architectonic structures inspired by natural shapes is now expanding (\cite{ch,he}), and the designs of singular minimal surfaces give stability in these constructions. Finally, in the future, it would be desirable to investigate the implementation of methods of discrete differential geometry which can produce this type of roof model in practice.

\section*{Acknowledgment.}
The author wishes to thank the  anonymous referees for a careful reading of the manuscript,  providing many useful remarks  and corrections. These suggestions greatly helped to improve the final version. In particular, one of the referees pointed out  the reference \cite{ni} of Nitsche for the Example \ref{ex1}.   The author also thanks to Bennet Palmer, Alvaro P\'ampano and Anthony Gruber who revised the initial draft. This work has been partially supported by the grant no. PID2020-117868GB-I00 Ministerio de Ciencia e Innovaci\'on.
  

  
 Rafael L\'opez works in classical differential geometry, in particular, surfaces with prescribed mean curvature. He is a Professor of Mathematics at the University of Granada where he received his Ph.D. in 1996. Rafael  enjoys performing mathematics outreach activities in schools using soap bubbles and, in his spare time, he likes  trekking and cycling in Sierra Nevada.

\end{document}